\pdfoutput=1
\RequirePackage{ifpdf}
\ifpdf 
\documentclass[pdftex]{sigma}
\else
\documentclass{sigma}
\fi

\DeclareMathOperator\Ai{{Ai}}
\DeclareMathOperator\Bi{{Bi}}

\numberwithin{equation}{section}

\newtheorem{Theorem}{Theorem}[section]
\newtheorem{Lemma}[Theorem]{Lemma}
\newtheorem{Proposition}[Theorem]{Proposition}
 { \theoremstyle{definition}
\newtheorem{Remark}[Theorem]{Remark} }

\begin{document}

\newcommand{\arXivNumber}{1408.5654}

\allowdisplaybreaks

\renewcommand{\PaperNumber}{070}

\FirstPageHeading

\ShortArticleName{Uniform Asymptotics of Orthogonal Polynomials Arising from Coherent States}

\ArticleName{Uniform Asymptotics of Orthogonal Polynomials\\ Arising from Coherent States}

\Author{Dan DAI~$^\dag$, Weiying HU~$^\dag$ and Xiang-Sheng WANG~$^\ddag$}

\AuthorNameForHeading{D.~Dai, W.~Hu and X.-S.~Wang}

\Address{$^\dag$~Department of Mathematics, City University of Hong Kong, Hong Kong}
\EmailD{\href{mailto:dandai@cityu.edu.hk}{dandai@cityu.edu.hk}, \href{mailto:weiyinghu2-c@my.cityu.edu.hk}{weiyinghu2-c@my.cityu.edu.hk}}

\Address{$^\ddag$~Department of Mathematics, Southeast Missouri State University,\\
\hphantom{$^\ddag$}~Cape Girardeau, MO 63701, USA}
\EmailD{\href{mailto:xswang@semo.edu}{xswang@semo.edu}}
\URLaddressD{\url{http://faculty.semo.edu/xswang/}}

\ArticleDates{Received April 01, 2015, in f\/inal form August 25, 2015; Published online August 31, 2015}

\Abstract{In this paper, we study a family of orthogonal polynomials $\{\phi_n(z)\}$ arising from nonlinear coherent states in quantum optics. Based on the three-term recurrence relation only, we obtain a uniform asymptotic expansion of $\phi_n(z)$ as the polynomial degree $n$ tends to inf\/inity. Our asymptotic results suggest that the weight function associated with the polynomials has an unusual singularity, which has never appeared for orthogonal polynomials in the Askey scheme. Our main technique is the Wang and Wong's dif\/ference equation method. In addition, the limiting zero distribution of the polynomials $\phi_n(z)$ is provided.}


\Classification{41A60; 33C45}

\section{Introduction}

Coherent states were f\/irst studied by Schr\"{o}dinger \cite{Schr} in the early years of quantum mechanics. In 1960s,
they were rediscovered by Glauber \cite{Glauber1,Glauber2}, Klauder \cite{Klauder1,Klauder2} and Sudarshan \cite{Sudarshan}
in the study of quantum optics. Since then, coherent states and their generalizations have been used in nearly all branches of quantum physics,
for example, nuclear, atomic and condensed matter physics, quantum optics, quantum f\/ield theory, quantization and dequantization problems, etc.
For more properties and applications of coherent states, one may refer to \cite{Ali:Ant:Gaz,Com:Rob,Gazeau:book}.

Recently, Ali and Ismail \cite{Ali:Ism} studied two sets of orthogonal polynomials associated with a~fa\-mi\-ly of nonlinear coherent states in quantum optics. To construct the f\/irst set of orthogonal polynomials, they derived a measure from resolution of the identity for the coherent states and demonstrated that this measure, denoted by $d\mu$, can be extended to an even probability measure supported on a symmetric interval $[-M,M]$ with $0< M \leq \infty$; see \cite[equation~(1.14)]{Ali:Ism}. The moments of $d\mu$ are def\/ined as
\begin{gather*}
  \mu_n:=\int_{-M}^M x^n d\mu(x), \qquad n = 0, 1, 2, \dots.
\end{gather*}
By symmetry, it is easily seen that $\mu_{2n+1} = 0$. The f\/irst set of monic orthogonal polynomials $\{ \psi_n(x)\colon \psi_n(x)= x^n + \cdots \}$ is constructed as
\begin{gather*}
  \psi_n(x):= \frac{1}{D_{n-1}} \left| \begin{matrix}
    \mu_0 & \mu_1 & \cdots & \mu_n \\
    \vdots & \vdots & \cdots & \vdots \\
    \mu_{n-1} & \mu_{n} & \cdots & \mu_{2n-1} \\
    1 & x & \cdots & x^n
  \end{matrix} \right|,
\end{gather*}
where $D_n$ is the Hankel determinant
\begin{gather*}
  D_n:=  \det \left( \mu_{j+k} \right)_{j,k=0}^{n}.
\end{gather*}
It is easily seen that $\psi_n(x)$'s are orthogonal with respect to the measure $d\mu$ on $[-M,M]$.

To def\/ine the other set of orthogonal polynomials, Ali and Ismail \cite{Ali:Ism} introduced the ratio of moments
\begin{gather} \label{an-mun}
  \lambda_n = \frac{\mu_{2n}}{\mu_{2n-2}},
\end{gather}
and showed that the sequence $\{ \lambda_n\}$
is strictly increasing. Furthermore, $\{ \lambda_n\}$ is bounded by $M^2$ if $M<\infty$; and unbounded if $M =\infty$.
The second set of orthogonal polynomials $\{\phi_n(x)\}$
is generated from the following three-term recurrence relation
\begin{gather} \label{phi-n-def}
  x \phi_n(x) = \sqrt{\frac{\lambda_{n+1}}{2}} \phi_{n+1}(x) + \sqrt{\frac{\lambda_n}{2}} \phi_{n-1}(x), \qquad n = 1,2, \dots,
\end{gather}
with initial conditions $\phi_0(x) =1 $ and $\phi_1(x) = \sqrt{\frac{2}{\lambda_1}}  x$. According to the Favard's theorem
(see, e.g., Szeg\H{o}~\cite{Sze}), there exists an even probability measure~$d\mu^*$
supported on the real line~$\mathbb{R}$ such that~$\phi_n(x)$'s are the corresponding orthonormal polynomials. Ali and Ismail~\cite{Ali:Ism} also pointed out that, in general, the measure~$d\mu^*$ is dif\/ferent from the orthogonality measure $d\mu$ for the f\/irst set of polynomials~$\psi_n(x)$.  They studied several examples in~\cite{Ali:Ism}, but no obvious relations between~$d\mu$ and~$d\mu^*$ have been concluded. It would be interesting to f\/ind an explicit relation between these two measures, as well as the corresponding two sets of polynomials. As the second measure~$d\mu^*$ can not be explicitly given in many general cases, one can only obtain some information (such as asymptotic behaviors) of~$\phi_n(x)$ directly from their three-term recurrence relation~\eqref{phi-n-def}.

Let us consider an example from Ali and Ismail \cite{Ali:Ism}. Let $\alpha$ and $\beta$ be two nonnegative constants such that
\begin{gather*}
  0\le\alpha<\beta,
\end{gather*}
and $w_\psi(x)$ be the weight function for the f\/irst set of polynomials $\psi_n(x):=\psi_n(x;\alpha,\beta)$
\begin{gather} \label{psi-weight}
  w_\psi(x):= w_\psi(x,\alpha,\beta) = \frac{2^\beta \Gamma(\beta + \frac{1}{2}) }{ \sqrt{\pi}   \Gamma(\beta+\alpha) \Gamma(\beta - \alpha)} e^{-x^2} K_\alpha\big(x^2\big) |x|^{2\beta-1},  \qquad x\in \mathbb{R},
\end{gather}
where $K_\alpha(x)$ is the modif\/ied Bessel function. The moments of the measure $w_\psi(x)dx$ can be computed explicitly as follows
\begin{gather*}
  \mu_{2n} = \frac{(\beta+\alpha)_n (\beta-\alpha)_n} {2^n (\beta+\frac{1}{2})_n},
\end{gather*}
where $(a)_n := a (a+1) \cdots (a+n-1)$ is  the Pochhammer symbol. According to~\eqref{an-mun}, the recurrence coef\/f\/icient~$\lambda_n$ for the second set of orthogonal polynomials $\phi_n(x):=\phi_n(x;\alpha,\beta)$ in~\eqref{phi-n-def} is given by
\begin{gather} \label{an-def}
  \lambda_n: = \lambda_n(\alpha,\beta) = \frac{(n+ \beta + \alpha -1)(n+ \beta - \alpha  -1)}{2 (n+ \beta -\frac{1}{2})}, \qquad n \geq 1;
\end{gather}
see \cite[equations (3.13) and (3.14)]{Ali:Ism}. The recurrence coef\/f\/icient $\lambda_n$ is a rational function in $n$ when $\alpha \neq \frac{1}{2}$, and can be rewritten as
\begin{gather} \label{an-new-expression}
   \lambda_n = \frac{1}{2}\left(n+\beta-\frac{3}{2} - \frac{\alpha^2-\frac{1}{4}}{n+\beta-\frac{1}{2}} \right).
\end{gather}
Clearly from the above expression, when $\beta=\frac{3}{2}$ and $\alpha = \frac{1}{2}$, $\lambda_n$ reduces to a monomial, i.e., $\lambda_n = \frac{n}{2}$. In this case, $\phi_n(x)$'s are indeed the Hermite polynomials
\begin{gather} \label{phi-n-hermite}
  \phi_n(x) = \frac{1}{\sqrt{2^n    n!}}   H_n\big(\sqrt{2}   x\big),
\end{gather}
and their weight function is $\sqrt{\frac{2}{\pi}}   e^{-2x^2}$, $x \in \mathbb{R}$. Therefore, for general~$\alpha$ and~$\beta$, the polynomials~$\phi_n(x)$ can be viewed as a perturbation of the Hermite polynomials.

\begin{Remark}
  As mentioned in the previous paragraph, the relations between the two sets of polynomials $\psi_n(x)$ and $\phi_n(x)$ are unclear for general $\alpha$ and $\beta$. It is interesting to see that, when $\beta=\frac{3}{2}$ and $\alpha = \frac{1}{2}$, the weight functions $w_\psi(x)$ and $w_\phi(x)$ have the following simple relation
  \begin{gather*}
    w_\psi(x) = 2 |x|  e^{-2x^2} = \sqrt{2\pi}  |x|w_\phi(x) , \qquad x \in \mathbb{R}.
  \end{gather*}
  Here we make use of~\eqref{psi-weight} and the fact that $K_{\frac{1}{2}}(z) = \sqrt{\frac{\pi}{2z}} e^{-z}$; see \cite[equation~(10.39.2)]{DLMF}.
  For the case $\alpha=\frac{1}{2}$ and $\beta>\alpha$ is arbitrary, $\phi_n(x)$ is the associated Hermite polynomials introduced by Askey and Wimp~\cite{Ask:Wimp}.
  One can show that
  \begin{gather*}
    w_\psi(x)=\frac{2^{\beta-1/2}}{\Gamma(\beta-1/2)} |x|^{2\beta - 2} e^{-2 x^2} ,\qquad
    w_\phi(x)=\frac{\sqrt{2/\pi}}{\Gamma(\beta-1/2)}|U(\beta-2,2ix)|^{-2},
  \end{gather*}
  where $U(a,x)$ is the parabolic cylinder function.
  In view of the large-variable asymptotic formula for parabolic cylinder function \cite[equation~(12.9.1)]{DLMF}, we have as $x\to\infty$,
  \begin{gather*}
    w_\phi(x)\sim \frac{\sqrt{2/\pi}}{\Gamma(\beta-1/2)} |2x|^{2\beta-3} e^{-2x^2}=\frac{2^{\beta-2}}{\sqrt{\pi}   |x|}w_\psi(x).
  \end{gather*}
\end{Remark}

In this paper, we focus on asymptotics of the second set of orthogonal polynomials $\phi_n(x)$ as their degree $n$ tends to inf\/inity. According to Ali and Ismail \cite{Ali:Ism}, nothing is known about the polynomials~$\phi_n(x)$ except for the
three-term recurrence relation. Therefore, some existing asymptotic methods, such as dif\/ferential equation methods in Olver~\cite{olver}, integral methods in Wong~\cite{Wong} and Riemann--Hilbert methods in Deift~\cite{deift}, cannot be applied. As a consequence, one can only make use of the three-term recurrence relation~\eqref{phi-n-def} to study the asymptotic properties of~$\phi_n(x)$.

In the literature, the asymptotic study of orthogonal polynomials via the three-term recurrence relation has been intensively investigated. For example, when the recurrence coef\/f\/icients converge to bounded constants, M\'{a}t\'{e} et al.~\cite{Mat:Nev:Tot}, Van Assche and Geronimo~\cite{VanAssche:Geronimo1988} obtained the
asymptotics of the orthogonal polynomials on and of\/f the essential spectrum of the orthogonality measure. When the recurrence coef\/f\/icients are regularly and slowly varying, the asymptotics away from the oscillatory were given by Van Assche and Geronimo~\cite{VanAssche:Geronimo1989} and Geronimo et al.~\cite{Ger:Smi:VAs}. Recently,
asymptotic techniques are also extended to study higher-order three-term recurrence relations which are satisf\/ied by some multiple orthogonal polynomials; see Aptekarev et al.~\cite{Apt:Kal:Saff}. Note that, as all results mentioned here are only valid either in the exponential or oscillatory regions, they are not uniformly valid in the neighborhood of the smallest or largest zeros (namely, transition points) of the polynomials.

In recent years, a turning point theory for second-order dif\/ference equations has been developed by Wong and his colleagues in~\cite{Cao:Li,wang-wong2003,wang-wong2005}. Their theory can be viewed as an analogue of the asymptotic theory for linear second-order dif\/ferential equations; see the def\/initive book by Olver \cite{olver}. In their papers~\cite{Cao:Li,wang-wong2003,wang-wong2005}, two linearly independent solutions are derived, which are uniformly valid in the neighborhood of the transition points. Depending on properties of the transition points, Airy-type or Bessel-type asymptotic expansions emerge. Since one can treat the recurrence relation \eqref{phi-n-def} as a second-order linear dif\/ference equation, the asymptotics of the polynomials can be obtained by using Wang and Wong's method if one can determine the coef\/f\/icients for the two linearly independent solutions. This is done via the method in Wang and Wong~\cite{xswang2012}, which makes use of the recurrence relation only; see also~\cite{xswang2014}. Recently, Dai, Ismail and Wang~\cite{Dai:Ism:Wang} successfully apply the dif\/ference equation method to derive uniform asymptotic expansions for some orthogonal polynomials from indeterminate moment problems. It should be mentioned that Geronimo~\cite{Geronimo} also studied WKB approximations for dif\/ference equations and obtained the Airy-type expansions.

Based on Wang--Wong's dif\/ference equation method~\cite{wang-wong2003,wang-wong2005} and the matching technique in~\cite{xswang2012}, we will derive a uniform asymptotic formula for the orthogonal polynomials $\phi_n(x)$ from their recurrence relation. Furthermore, we adopt the method of Kuijlaars and Van Assche~\cite{Kui:Ass} to f\/ind the limiting zero distribution of~$\phi_n(x)$ from the recurrence relation, too. Note that Coussement et al.~\cite{Cou:Cou:VAs} recently extended the method in Kuijlaars and Van Assche~\cite{Kui:Ass} to derive the limiting zero distribution for polynomials generated by a four-term recurrence relation. In this case, the polynomials are a class of multiple orthogonal polynomials.

There are mainly two reasons which motivate our research in this paper. First, as nothing is known about the polynomials~$\phi_n(x)$ except for the three-term recurrence relation, obtaining asymptotic behaviors and limiting zero distributions of~$\phi_n(x)$ as~$n\to\infty$ are important steps for us to understand the properties of such kind of orthogonal polynomials arising from coherent states. Second, we notice that the term $\frac{\alpha^2-\frac{1}{4}}{n+\beta-\frac{1}{2}}$ in~\eqref{an-new-expression} seems rare because such kind of~$\frac{1}{n}$-term never appears in the recurrence coef\/f\/icients for the hypergeometric polynomials in the Askey scheme. A further investigation shows that this term will contribute for a power of $x^{1/x^2}$ in the asymptotic formula of the polynomials. Although an explicit formula of the weight function for~$\phi_n(x)$ can not be determined at this stage, asymptotic results suggest that this~$\frac{1}{n}$-term in the coef\/f\/icients of recurrence relation may induce a $|x|^{1/x^2}$-term in the weight function for the corresponding orthogonal polynomials. Note that this phenomenon does not appear for any orthogonal polynomials in the Askey scheme. More detailed discussions can be found in Section~\ref{sec:discuss}.

The rest of the paper is arranged as follows. In Section~\ref{sec:results}, we state our main results with uniform asymptotic formula and limiting zero distribution for the orthogonal polynomials~$\phi_n(x)$. In Section~\ref{sec:wang-wong}, we introduce Wang and Wong's dif\/ference equation method. In Section~\ref{sec:proof}, we derive some ratio and non-uniform asymptotic formulas for~$\phi_n(x)$. The proof of our main theorem is also given in this section. In Section~\ref{sec:discuss}, we conclude our paper with some discussions together with suggestions about several possible problems for the future work.

\section{Main results} \label{sec:results}

In this paper, we will follow the general framework developed in \cite{Dai:Ism:Wang} to derive uniform asymptotic formulas.
Before stating our main theorem, we shall introduce some constants and functions. Def\/ine
\begin{gather} \label{kn-def}
  k_n := k_n(\alpha,\beta) = \sqrt{ \frac{\Gamma\left( \frac{n+\beta + \alpha}{2} \right) \Gamma\left( \frac{n+\beta - \alpha}{2} \right) \Gamma\left( \frac{n+\beta + 3/2}{2} \right)}
  {\Gamma\left( \frac{n+\beta + \alpha+1}{2} \right) \Gamma\left( \frac{n+\beta - \alpha + 1}{2} \right) \Gamma\left( \frac{n+\beta + 1/2}{2} \right)} }
\end{gather}
and
\begin{gather}
  \frac{2}{3} [U(t)]^{\frac{3}{2}}  :=  t\sqrt{t^2 -1 } - \log\big(t + \sqrt{t^2 -1 }\big), \qquad t \geq 1, \label{Ut-def1} \\
  \frac{2}{3} [-U(t)]^{\frac{3}{2}}  :=   \cos^{-1} t - t \sqrt{1-t^2}, \qquad   -1<t < 1. \label{Ut-def2}
\end{gather}
\begin{Remark}
  The function $U(t)$ def\/ined above is analytic for $t$ in a complex neighborhood of~1. Moreover, we have the following asymptotic formula
   \begin{gather*}
     U(t) = 2(t-1) + O(t-1)^2, \qquad \textrm{as } t \to 1.
   \end{gather*}
   Especially, $U(t)$ is monotonically increasing for $t\in[1-\delta,1+\delta]$ with $\delta>0$ being a small positive number.
\end{Remark}

One can easily see from \eqref{phi-n-def} that $\phi_{n}(x)$ is symmetric with respect to $x$, i.e., $\phi_{2n}(x)$ and $\phi_{2n+1}(x)$ are even and odd functions, respectively. Thus, it suf\/f\/ices to investigate~$\phi_n(x)$ for positive~$x$ only. Our f\/irst main result is stated in the following theorem and its proof will be postponed until the end of Section~\ref{sec:proof}.
\begin{Theorem} \label{main-thm}
  Let $N = n+\beta-1$. With $k_n$ and $U(t)$ defined in \eqref{kn-def}, \eqref{Ut-def1} and \eqref{Ut-def2}, respectively, we have as $n\to\infty$,
  \begin{gather}
    \phi_{n}\big(N^{\frac{1}{2}}   t\big)   \sim
     k_n \pi^{\frac{1}{4}} \sqrt{ \frac{\Gamma(\beta+\alpha) \Gamma(\beta-\alpha)}{ \Gamma(\beta+\frac{1}{2}) } } \frac{e^{N t^2  + \frac{(\alpha^2-\frac{1}{4}) \log (Nt^2) }{4Nt^2}  } }{ \big( 2 N^{\frac{1}{2}}   t\big) ^{\beta - \frac{3}{2}} } \left(\frac{ U(t)}{t^2 - 1}\right)^{\frac{1}{4}} \nonumber \\
\hphantom{\phi_{n}\big(N^{\frac{1}{2}}   t\big)   \sim }{}
       \times \left[ \Ai\big(N^{\frac{2}{3}} U(t)\big) \sum_{s=0}^\infty \frac{\tilde A_s(U)}{N^{s-\frac{1}{6}}} + \Ai'\big(N^{\frac{2}{3}} U(t)\big) \sum_{s=0}^\infty \frac{\tilde B_s(U)}{N^{s+\frac{1}{6}}} \right], \label{main-for}
  \end{gather}
  uniformly for $t$ in $[\delta,\infty)$, where $\delta>0$ is any fixed small positive number. The leading coefficients $\tilde{A}_0(U)$ and $\tilde{B}_0(U)$ are given by
\begin{gather*}
  \tilde{A}_0(U) = 1, \qquad \tilde{B}_0(U) = 0,
\end{gather*}
and the higher coefficients $\tilde{A}_s(U)$ and $\tilde{B}_s(U)$ can be determined recursively.
The asymptotic formula for $t\in(-\infty,-\delta]$ is easily obtained via the reflection formula $\phi_n(-x)=(-1)^n\phi_n(x)$.
\end{Theorem}

\begin{Remark}
  When $\beta=\frac{3}{2}$ and $\alpha = \frac{1}{2}$, our expansion \eqref{main-for} reduces to
  \begin{gather*}
    \phi_{n}\big(N^{\frac{1}{2}}   t\big)   \sim   \pi^{\frac{1}{4}} \sqrt{ \frac{\Gamma\big(\frac{n+1}{2}\big) }{ \Gamma\big(\frac{n+2}{2}\big) } } e^{N t^2}  \left(\frac{ U(t)}{t^2 - 1}\right)^{\frac{1}{4}}  \nonumber \\
    \hphantom{\phi_{n}\big(N^{\frac{1}{2}}   t\big)   \sim}{}
 \times \biggl[ \Ai\big(N^{\frac{2}{3}} U(t)\big) \sum_{s=0}^\infty \frac{\tilde A_s(U)}{N^{s-\frac{1}{6}}} + \Ai'\big(N^{\frac{2}{3}} U(t)\big) \sum_{s=0}^\infty \frac{\tilde B_s(U)}{N^{s+\frac{1}{6}}} \biggr]
  \end{gather*}
    with $N = n + \frac{1}{2}$. Moreover, one can show that $\tilde{A}_{2s+1}(U) = \tilde{B}_{2s}(U) = 0$ for $\alpha=1/2$ and $\beta=3/2$. On account of \eqref{phi-n-hermite}, the above formula agrees with the asymptotic expansion for the Hermite polynomials in \cite[equation~(12.10.35)]{DLMF}.
  For the general values of~$\alpha$ and~$\beta$, the coef\/f\/icients~$\tilde A_s(U)$ and~$\tilde B_s(U)$ can be obtained explicitly using the recursive formulas (4.34)--(4.37) in~\cite{wang-wong2003}. But, the resulting formulas for the higher-order coef\/f\/icients (that is, $\tilde{A}_{s}(U)$ and $\tilde{B}_{s}(U)$ with $s\ge1$) would be too complicated and we omit the details here. The readers are referred to \cite[Section~4]{wang-wong2003} for more information about the recursive formulas for the coef\/f\/icients. We remark that the leading term of the asymptotic formula \eqref{main-for} is simple since $\tilde A_0(U)=1$ and $\tilde B_0(U)=0$.
\end{Remark}

From the three-term recurrence relation, we can also obtain asymptotic zero distribution of the rescaled polynomials $\phi_{n}(m^{\frac{1}{2}}  t)$; as stated in the following theorem.

\begin{Theorem}
    Let $\nu(p_n)$ be the normalized zero counting measure
  \begin{gather*}
    \nu(p_n):=\frac{1}{n} \sum_{j=1}^n \delta_{x_{j,n}},
  \end{gather*}
  where $x_{j,n}$, $j=1,\dots, n$ are zeros of $p_n(x)$ and $\delta_{x_{j,n}}$ denotes the Dirac point mass at the zero~$x_{j,n}$. Let $\Phi_{n,m}(t)$ be the rescaled polynomials
  \begin{gather*}
    \Phi_{n,m}(t) : = \phi_n\big(\sqrt{m} t\big).
  \end{gather*}
  Then, for every $c>0$,
  \begin{gather*}
    \lim_{{n/m\to c,}\atop n \to \infty} \nu(\Phi_{n,m}) = \nu_c
  \end{gather*}
  exists and has the density
  \begin{gather} \label{limiting-zero}
    \frac{d\nu_c }{dt} = \begin{cases}
      \displaystyle\frac{2}{\pi c} \sqrt{c-t^2}, & t \in \big({-}\sqrt{c},\sqrt{c}\big), \\
      0, & \textrm{elsewhere}.
    \end{cases}
  \end{gather}
\end{Theorem}

\begin{proof}
  The proof is an easy application of Theorem~1.4 in Kuijlaars and Van Assche~\cite{Kui:Ass}. From~\eqref{phi-n-def}, we have
  \begin{gather*}
    t \Phi_{n,m}(t) = \sqrt{\frac{\lambda_{n+1}}{2m}} \Phi_{n+1,m}(t) + \sqrt{\frac{\lambda_{n+1}}{2m}} \Phi_{n-1,m}(t).
  \end{gather*}
  Then, when $n/m \to c$, the quantities in \cite[equations~(1.7) and (1.8)]{Kui:Ass} are
  \begin{gather*}
    a(c) = \frac{\sqrt{c}}{2}, \qquad b(c) = 0, \qquad \alpha(c) = -\sqrt{c}, \qquad \beta(c) = \sqrt{c}.
  \end{gather*}
  The density in \eqref{limiting-zero} immediately follows from \cite[equation~(1.9)]{Kui:Ass}, i.e.,
  \begin{gather*}
    \frac{d\nu_c }{dt} = \frac{1}{c} \int_{t^2}^{c} \frac{1}{\pi \sqrt{(\beta(s)-t)(t-\alpha(s))}}ds = \frac{1}{\pi c} \int_{t^2}^{c} \frac{1}{\sqrt{s-t^2}}ds = \frac{2}{\pi c} \sqrt{c-t^2}.
  \end{gather*}
  This f\/inishes the proof of the theorem.
\end{proof}

\begin{Remark}
  One can see that the limiting zero distribution \eqref{limiting-zero} is independent of the parame\-ters~$\alpha$ and~$\beta$, and it coincides with that of the Hermite polynomials. This phenomenon again suggests that the polynomials~$\phi_n(x)$ can be viewed as a perturbation of the Hermite polynomials.
\end{Remark}

\section{Wang and Wong's dif\/ference equation method} \label{sec:wang-wong}

To apply Wang and Wong's dif\/ference equation method, we f\/irst let
\begin{gather} \label{phin-pn}
  \phi_{n}(x) = k_n \, p_n(x)
\end{gather}
and transform the recurrence relation \eqref{phi-n-def} into their standard form
\begin{gather} \label{standard-form}
  p_{n+1}(x) - (A_n x + B_n) \, p_n(x) + p_{n-1}(x) = 0.
\end{gather}
From the def\/inition of $k_n$ in \eqref{kn-def}, it is easily verif\/ied that $k_{n+1} = \sqrt{ \frac{\lambda_n}{\lambda_{n+1}}} k_{n-1} $.
Then, \eqref{phi-n-def} is reduced to the standard form \eqref{standard-form} with $B_n=0$ and
\begin{gather} \label{An-asy}
  A_n = \sqrt{ \frac{2}{\lambda_n} } \frac{k_n}{k_{n-1}}  \sim n^{-\theta} \sum_{s=0}^\infty \frac{a_s}{n^s}, \qquad  \theta=\frac{1}{2}, \qquad a_0 = 2, \qquad   a_1 =  1-\beta,
\end{gather}
as $ n \to \infty$.
Since $A_n$ is of order $O(n^{-\frac{1}{2}})$ when $n$ is large, to balance the term $A_n x $ in~\eqref{standard-form}, we introduce a new scale $x = N^{\frac{1}{2}} t$ with $N = n + \tau_0$, where $\tau_0$ is a constant to be determined. The characteristic equation for~\eqref{standard-form} is
\begin{gather*}
  \xi^2- a_0 t  \xi + 1= 0
\end{gather*}
with $a_0=2$; see \eqref{An-asy}. The roots of this equation are
\begin{gather} \label{roots-chara-eqn}
  \xi(t) = t \pm \sqrt{t^2 - 1},
\end{gather}
and they coincide when the quantity inside the above square root vanishes, that is, $t = t_\pm = \pm 1$. In general, the points $t_\pm$ separate the oscillatory interval from the exponential intervals, and the asymptotic behaviors of solutions to \eqref{standard-form} change dramatically when~$t$ passes through $t_\pm$. So, the points $t_\pm$ are called transition points for the dif\/ference equation~\eqref{standard-form} by Wang and Wong in~\cite{wang-wong2003, wang-wong2005}.

By symmetry, we only need to consider the right transition point $t_+=1$. According to the main theorem in \cite[p.~189]{wang-wong2003}, we have the following Airy-type expansion.

\begin{Proposition}\label{prop-uniform}
  When $n$ is large, $p_n(x)$ in \eqref{standard-form} can be expressed as
  \begin{gather*}
   p_n(x) = C_1(x) P_n(x) + C_2(x) Q_n(x),
  \end{gather*}
  where $C_1(x)$ and $C_2(x)$ are functions independent of $n$, and $P_n(x)$ and $Q_n(x)$ are two linearly independent solutions of \eqref{standard-form} satisfying the following Airy-type asymptotic expansions in the neighborhood of $t_+=1$
  \begin{gather*}
  P_n\big(N^{\frac{1}{2}} t\big) \sim  \left(\frac{U(t)}{t^2 - 1}\right)^{\frac{1}{4}} \left[ \Ai\big(N^{\frac{2}{3}} U(t)\big) \sum_{s=0}^\infty \frac{\tilde{A}_s(U)}{N^{s-\frac{1}{6}}} + \Ai'\big(N^{\frac{2}{3}} U(t)\big) \sum_{s=0}^\infty \frac{\tilde{B}_s(U)}{N^{s+\frac{1}{6}}} \right]
\end{gather*}
and
\begin{gather*}
  Q_n\big(N^{\frac{1}{2}} t\big) \sim \left(\frac{U(t)}{t^2 - 1}\right)^{\frac{1}{4}} \left[ \Bi\big(N^{\frac{2}{3}} U(t)\big) \sum_{s=0}^\infty \frac{\tilde{A}_s(U)}{N^{s-\frac{1}{6}}} + \Bi'\big(N^{\frac{2}{3}} U(t)\big) \sum_{s=0}^\infty \frac{\tilde{B}_s(U)}{N^{s+\frac{1}{6}}} \right],
\end{gather*}
where $N = n+ \beta -1$ and $U(t)$ is defined in \eqref{Ut-def1} and \eqref{Ut-def2}.
Here the leading coefficients are given by
\begin{gather*}
  \tilde{A}_0(U) = 1, \qquad \tilde{B}_0(U) = 0.
\end{gather*}
\end{Proposition}

\begin{proof}
   Recall the asymptotic expansion for $A_n$ in \eqref{An-asy}, the equation \eqref{standard-form} falls into the case $\theta \neq 0$ and $t_+ \neq 0 $ considered in~\cite{wang-wong2003}. Following their approach, we choose
\begin{gather*}
  \tau_0 = - \frac{a_1 t_+ }{2 \theta} = \beta - 1 \qquad \textrm{and} \qquad N = n + \tau_0.
\end{gather*}
Then our proposition follows from the main theorem in~\cite{wang-wong2003}.
\end{proof}


\section{Ratio and non-uniform asymptotics} \label{sec:proof}

To determine the coef\/f\/icients $C_1(x)$ and $C_2(x)$ in Proposition~\ref{prop-uniform}, more asymptotic information about the polynomials $\phi_n(x)$ is required. Note that, people had to obtain the asymptotic information by using results derived from other methods, see, e.g., \cite{wang-wong2003,wang-wong2005}. Here, through the approach developed by Wang and Wong in \cite{xswang2012}, we can get the non-uniform asymptotics we need. As a consequence, we obtain our main results in Theorem~\ref{main-thm} from the three-term recurrence relation~\eqref{phi-n-def} only.

Obviously, $\phi_n(x) = \gamma_n x^n + \cdots $ def\/ined in \eqref{phi-n-def} are orthonormal polynomials. It is easily seen from \eqref{phi-n-def} that the leading coef\/f\/icient of $\phi_n(x)$ is given by
\begin{gather*}
    \gamma_n := \sqrt{\frac{2^n}{\lambda_1 \cdots \lambda_n}} = 2^n \sqrt{\frac{\Gamma(\beta + \alpha) \Gamma(\beta -\alpha) \Gamma(n + \beta + \frac{1}{2})}{\Gamma(n + \beta + \alpha) \Gamma(n + \beta -\alpha) \Gamma(\beta + \frac{1}{2})}}.
\end{gather*}
For convenience, we f\/irst consider the corresponding monic polynomials $\pi_n(x)$, that is,
\begin{gather} \label{pin-phin}
  \pi_n(x) = x^n + \cdots = \frac{1}{\gamma_n} \phi_n(x).
\end{gather}
As the polynomials $\pi_n(x)$, $\phi_n(x)$ and $p_n(x)$ in~\eqref{standard-form} only dif\/fer by some constant factors, one can see that most of the zeros of $\pi_n(\sqrt{N} t)$ are asymptotically contained in the interval $[-1,1]$ from~\eqref{roots-chara-eqn}.
Then it is possible to consider the ratio $\frac{\pi_k(\sqrt{N} t)}{\pi_{k-1}(\sqrt{N} t)}$ outside the oscillatory inter\-val~$[-1,1]$. The ratio asymptotics for $\pi_k(x)$ is given as follows.

\begin{Lemma}
  Let $N = n+ \beta -1$ and $ w_k(x)$ be defined as
  \begin{gather} \label{wk-def}
    w_k(x):= \frac{\pi_k(x)}{\pi_{k-1}(x)}, \qquad k = 1, \dots , n.
  \end{gather}
  Then with $\lambda_n$ given in \eqref{an-def}, we have
  \begin{gather} \label{wk-asy}
    w_k\big(\sqrt{N}   z \big) = \sqrt{N} \frac{z + \sqrt{z^2 - \frac{2\lambda_k}{N}}}{2} \left[ 1+ \frac{\lambda_k - \lambda_{k-1}}{2[N z^2 - 2\lambda_k]} + O\big(N^{-2}\big) \right] \qquad \textrm{for} \quad |z| > 1,
  \end{gather}
  where the error term $O(N^{-2}) $ is uniformly valid for all $k = 1, \dots,n$. Here the principal branch of the square root is taken such that $\sqrt{z^2 - \frac{2\lambda_k}{N}} \sim z$ as $z \to \infty$.
\end{Lemma}

\begin{proof}
  From \eqref{phi-n-def} and \eqref{pin-phin}, it is easily seen that the monic polynomial $\pi_n(x)$ satisf\/ies the following three-term recurrence relation
  \begin{gather*}
    x \pi_n(x) = \pi_{n+1}(x) + \frac{\lambda_n}{2} \pi_{n-1}(x)
  \end{gather*}
  with $\pi_0(x) = 1$ and $\pi_1(x) = x$. From \eqref{wk-def}, we have
  \begin{gather} \label{pi-wk}
    \pi_n(x) = \prod_{k=1}^n w_k(x).
  \end{gather}
  It then follows that $w_k(x)$ satisf\/ies a two-term nonlinear recurrence relation{\samepage
  \begin{gather*}
    w_{k+1}(x) = x - \frac{\lambda_k}{2 w_k(x)}, \qquad k \geq 1
  \end{gather*}
  with $w_1(x) = x$.}

  Next, we prove \eqref{wk-asy} by mathematical induction. When $k=1$, we have
  \begin{gather*}
    w_1\big(\sqrt{N}   z\big) =  \sqrt{N}    \frac{z + \sqrt{z^2 - \frac{2\lambda_1}{N}}}{2}     \frac{2z}{z + \sqrt{z^2 - \frac{2\lambda_1}{N}}}.
  \end{gather*}
  Since $|z| > 1$, we expand the denominator in the second fraction and obtain
  \begin{gather*}
    \frac{2z}{z + \sqrt{z^2 - \frac{2\lambda_1}{N}}} = \frac{2z}{2z - \frac{\lambda_1}{Nz} + O(N^{-2})} = 1 + \frac{\lambda_1}{2Nz^2} + O\big(N^{-2}\big).
  \end{gather*}
  As $\lambda_1$ is a constant, the above formula gives \eqref{wk-asy} when $k=1$.

  Assume \eqref{wk-asy} is true for $k=m$. It then follows
  \begin{gather*}
    w_{m+1}\big(\sqrt{N}   z\big)  =  \sqrt{N}  z - \frac{\lambda_m}{2w_m\big(\sqrt{N}   z\big)} \\
  \hphantom{w_{m+1}\big(\sqrt{N}   z\big)}{}
      =    \sqrt{N}   z - \frac{\lambda_m}{\sqrt{N} \left(z + \sqrt{z^2 - \frac{2\lambda_m}{N}}\right) } \left( 1 - \frac{\lambda_m - \lambda_{m-1}}{2[N z^2 - 2\lambda_m]} + O\big(N^{-2}\big) \right) \\
 \hphantom{w_{m+1}\big(\sqrt{N}   z\big)}{}
     =   \sqrt{N}  z - \frac{\sqrt{N} \left(z - \sqrt{z^2 - \frac{2\lambda_m}{N}}\right) }{2} \left( 1 - \frac{\lambda_m - \lambda_{m-1}}{2[N z^2 - 2\lambda_m]} + O\big(N^{-2}\big) \right) \\
 \hphantom{w_{m+1}\big(\sqrt{N}   z\big)}{}
     =  \sqrt{N} \frac{ z + \sqrt{z^2 - \frac{2\lambda_m}{N}} }{2} + \frac{\sqrt{N}(\lambda_m - \lambda_{m-1})\left(z - \sqrt{z^2 - \frac{2\lambda_m}{N}}\right)}{4[N z^2 - 2\lambda_m]} + O\big(N^{-3/2}\big).
  \end{gather*}
  One can see that, except the dif\/ference in~$\lambda_{m+1}$ and~$\lambda_m$, the f\/irst term in the above formula is very close to the leading term in~\eqref{wk-asy}. To get the exact formula, let us rewrite the above formula as
  \begin{gather} \label{wm-asy1}
    w_{m+1}\big(\sqrt{N}   z\big) =  \sqrt{N} \frac{z + \sqrt{z^2 - \frac{2\lambda_{m+1}}{N}}}{2} \left[ 1+ \frac{\lambda_{m+1} - \lambda_{m}}{2[N z^2 - 2\lambda_{m+1}]} \right] + R(z),
  \end{gather}
  where the remainder $R(z)$ is
  \begin{gather}
    R(z)  =  \sqrt{N} \frac{\sqrt{z^2 - \frac{2\lambda_{m}}{N}}-\sqrt{z^2 - \frac{2\lambda_{m+1}}{N}}}{2}  + \frac{\sqrt{N}(\lambda_m - \lambda_{m-1})\left(z - \sqrt{z^2 - \frac{2\lambda_m}{N}}\right)}{4[N z^2 - 2\lambda_m]} \nonumber \\
\hphantom{R(z)  =}{} - \frac{\sqrt{N}(\lambda_{m+1} - \lambda_{m})\left(z + \sqrt{z^2 - \frac{2\lambda_{m+1}}{N}}\right)}{4[N z^2 - 2\lambda_{m+1}]} + O\big(N^{-3/2}\big).\label{remainder-R}
  \end{gather}
  Note that, when $|z| > 1$, the quantity $\frac{z + \sqrt{z^2 - \frac{2\lambda_{m+1}}{N}}}{2}$ in \eqref{wm-asy1} is nonzero and uniformed bounded for all $m = 0, \dots, n-1$. To prove \eqref{wk-asy} for $k=m+1$, then it is suf\/f\/icient to show that~$R(z)$ is of order $O(N^{-3/2})$ for all $m\leq n-1$.

  Denote $\sigma_m := \lambda_{m+1} - \lambda_m$. By \eqref{an-def}, it is easily seen that
  \begin{gather*}
    \sigma_m = \frac{1}{2} + \frac{\alpha^2-\frac{1}{4}}{2\big(m+\beta-\frac{1}{2}\big)\big(m+\beta+\frac{1}{2}\big)},
  \end{gather*}
  which is uniformly bounded for all $m$. Let us rewrite~\eqref{remainder-R} as
  \begin{gather*}
    R(z)  =  \frac{1}{\sqrt{N}} \frac{\sigma_m}{\sqrt{z^2 - \frac{2\lambda_{m}}{N}}+\sqrt{z^2 - \frac{2(\lambda_m + \sigma_m)}{N}}}  + \frac{\sigma_{m-1}z }{4\sqrt{N}\left(z^2 - \frac{2\lambda_m}{N}\right)} - \frac{\sigma_{m-1}}{4\sqrt{N}\sqrt{z^2 - \frac{2\lambda_m}{N}}} \nonumber \\
 \hphantom{R(z)  =}{}       - \frac{\sigma_m z }{4\sqrt{N}\left( z^2 - \frac{2(\lambda_m+ \sigma_m)}{N}\right)} - \frac{\sigma_{m}}{4\sqrt{N}\sqrt{z^2 - \frac{2(\lambda_m+ \sigma_m)}{N}}}  + O\big(N^{-3/2}\big).
  \end{gather*}
  As $\frac{\sigma_m}{N} = O(N^{-1})$ uniformly for all $m$, the above formula gives
  \begin{gather*}
    R(z)  =  \frac{1}{\sqrt{N}} \frac{\sigma_m}{2\sqrt{z^2 - \frac{2\lambda_{m}}{N}} +O(N^{-1}) }  + \frac{\sigma_{m-1}z }{4\sqrt{N}\left(z^2 - \frac{2\lambda_m}{N}\right)} - \frac{\sigma_{m-1}}{4\sqrt{N}\sqrt{z^2 - \frac{2\lambda_m}{N}}} \nonumber \\
\hphantom{R(z)  =}{}  - \frac{\sigma_m z }{4\sqrt{N}\left( z^2 - \frac{2\lambda_m}{N} +O(N^{-1})\right)} - \frac{\sigma_{m}}{4\sqrt{N}\left( \sqrt{z^2 - \frac{2\lambda_m}{N}} +O(N^{-1})\right)}  + O\big(N^{-3/2}\big) \nonumber \\
\hphantom{R(z)  =}{}
     =   \frac{(\sigma_m-\sigma_{m-1}) \left(\sqrt{z^2 - \frac{2\lambda_m}{N}} - z\right)}{4\sqrt{N} \left(z^2 - \frac{2\lambda_m}{N} \right)} + O\big(N^{-3/2}\big).
  \end{gather*}
  When $m$ is bounded, $\lambda_m$ is bounded and $\sqrt{z^2 - \frac{2\lambda_m}{N}} - z = O(N^{-1})$, which means $R(z) = O(N^{-3/2})$. When $m=O(N^\delta), 0<\delta \leq 1$ is unbounded, we have
  \begin{gather*}
    \sqrt{z^2 - \frac{2\lambda_m}{N}} - z = O\big(N^{\delta-1}\big) \qquad \textrm{and} \qquad \sigma_m-\sigma_{m-1} = O\big(m^{-2}\big) = O\big(N^{-2\delta}\big).
  \end{gather*}
  Then, $R(z)$ is again of order $O(N^{-3/2})$. Thus, \eqref{wk-asy} holds for $k=m+1$ uniformly for all $m = 0, \dots, n-1$.

  This proves our lemma.
\end{proof}

The result in the above lemma is remarkable in the sense that \eqref{wk-asy} is valid not only when~$k$ is large but also when $k = 1, 2, \dots$. With this result, and taking~\eqref{pi-wk} into consideration, one can apply the trapezoidal rule or the Euler--Maclaurin formula to obtain the asymptotics of $\pi_n(\sqrt{N}  z)$ as $n \to \infty$, uniformly for~$z$ bounded away from~$[-1,1]$.

Before deriving the asymptotics of $\pi_n(\sqrt{N}   z)$, let us f\/irst compute the following summation for later use.

\begin{Lemma}\label{lemma-logn/n}
  For $\beta>0$, we have, as $n \to \infty$,
  \begin{gather} \label{sum-log-n}
    \sum_{k=1}^n \frac{1}{(k+\beta-\frac{1}{2})  \big(x^2 - k + x \sqrt{ x^2 - k} \big)}  =\frac{\log x^2}{ {2x^2}} + {O}\big(n^{-1}\big)
  \end{gather}
  for any complex $x$ such that $x/\sqrt n$ is bounded away from~$[-1,1]$. Here the principal branch of the square root is taken such that $\sqrt{x^2-k}\sim x$ as $x\to \infty$.
\end{Lemma}

\begin{proof}
Recall the Euler--Maclaurin formula in \cite[Theorem~6]{Wong}, we have
\begin{gather*}
  \sum_{k=1}^n {f(k)}= \int_{1}^n {f(t)}dt + \frac{1}{2}[f(n)+f(1)] + \sum_{s=1}^{m}\frac{B_{2s}}{(2s)!}\big[f^{(2s-1)}(n)-f^{(2s-1)}(1)\big] + R_m,
\end{gather*}
where $B_n$ and $B_n(t)$ are Bernoulli numbers and Bernoulli polynomials, respectively. The remainder $R_m$ is given explicitly below
\begin{gather*}
  R_m = \int_{1}^n \frac{B_{2m} - B_{2m} (t - [t])}{(2m)!}f^{(2m)}(t)dt
\end{gather*}
and satisf\/ies the following inequality
\begin{gather*}
  \left|R_m\right|\leqslant\big(2-2^{1-2m}\big)\frac{B_{2m}}{(2m)!}\int_{1}^n \big|f^{(2m)}(t)\big|dt.
\end{gather*}
To prove our formula \eqref{sum-log-n}, we def\/ine
\begin{gather} \label{f-def}
  f(t):=\frac{1}{\big(t+\beta-\frac{1}{2}\big)\big(x^2 + x\sqrt{x^2 - t}-t\big)}, \qquad 1 \leq t \leq n,
\end{gather}
and choose $m = 1$. Then the Euler--Maclaurin formula reduces to
\begin{gather} \label{EM-sum-m=1}
  \sum_{k=1}^n {f(k)} = \int_{1}^n {f(t)}dt + \frac{1}{2}[f(n)+f(1)] + \frac{B_2}{2} \big[f^{(1)}(n) - f^{(1)}(1)\big] + R_1.
\end{gather}
When $\beta\neq \frac12$, the f\/irst term in the above formula is calculated explicitly
\begin{gather*}
  \int_{1}^n {f(t)}dt  = \frac{1}{2\beta -1}\left[\frac{4x\operatorname{arctanh}{\left(\frac{\sqrt{2x^2 - 2n}}{\sqrt{2x^2 + 2\beta - 1}}\right)}}{\sqrt{\frac{1}{2}(2x^2 + 2\beta -1)}}- 2\log\frac{\big(x+ \sqrt{x^2 -n}\big)^2}{2n + 2\beta -1} \right]\nonumber \\
\hphantom{\int_{1}^n {f(t)}dt  =}{} - \frac{1}{2\beta-1}\left[\frac{4x\operatorname{arctanh}\left(\frac{\sqrt{2x^2 - 2}}{\sqrt{2x^2 + 2\beta - 1}}\right)}{\sqrt{\frac{1}{2}(2x^2 + 2\beta -1)}}- 2\log\frac{\big(x+ \sqrt{x^2 - 1}\big)^2}{2 + 2\beta -1} \right].
\end{gather*}
As $x/\sqrt n$ is bounded away from $[-1,1]$, then $x \to \infty$ when $n \to \infty$. Using the formula $\operatorname{arctanh}{z}=\frac12\log{\frac{1+z}{1-z}}$, we have
\begin{gather}
  \int_{1}^n {f(t)}dt = -\frac{1}{2x^2}\left(2\log\frac{x+\sqrt{x^2 - n}} {x+\sqrt{x^2 - 1}}-\log\frac{2n+2\beta-1} {2\beta+1} \right) +O\big(x^{-2}\big)\nonumber \\
 \hphantom{\int_{1}^n {f(t)}dt}{}
    =  \frac{\log x^2}{2x^2}-\frac{1}{x^2}\log\left(\frac{x}{\sqrt n}+\sqrt{\frac{x^2}{n} - 1}\right)+O\big(x^{-2}\big).\label{int-beta!=1/2}
\end{gather}
When $\beta=\frac12$, the f\/irst term in \eqref{EM-sum-m=1} can also be computed explicitly
\begin{gather}
  \int_{1}^n {f(t)}dt  = \frac{1}{x^2+x\sqrt{x^2-n}}+\frac{\log\big(x-\sqrt{x^2-n}\big)-\log(x+\sqrt{x^2-n})}{2x^2}\nonumber \\
\hphantom{\int_{1}^n {f(t)}dt  = }{} - \left[\frac{1}{x^2+x\sqrt{x^2-1}}+\frac{\log\big(x-\sqrt{x^2-1}\big)-\log\big(x+\sqrt{x^2-1}\big)}{2x^2}\right]\nonumber\\
\hphantom{\int_{1}^n {f(t)}dt }{}  = \frac{\log x^2}{2x^2}-\frac{1}{x^2}\log\left(\frac{x}{\sqrt n}+\sqrt{\frac{x^2}{n} - 1}\right)+O\big(x^{-2}\big).\label{int-beta=1/2}
\end{gather}
As $|x| > \sqrt{n}$, we have from \eqref{int-beta!=1/2} and \eqref{int-beta=1/2}
\begin{gather}\label{EM-leading}
  \int_{1}^n {f(t)}dt =\frac{\log{x^2}}{2x^2}+O\big(n^{-1}\big)
\end{gather}
uniformly for all $\beta>0$.

From the def\/inition of $f(t)$ in \eqref{f-def}, it is easily seen that
\begin{gather*}
  f^{(1)}(t)  =  \frac{-1}{\big(x^2 + x \sqrt{x^2 - t} - t\big) \big(t + \beta-\frac12 \big)^2} +\frac{1 + \frac{x}{2 \sqrt{x^2 - t}}}{\big(x^2 + x \sqrt{x^2 - t} - t\big)^2 \big( t + \beta-\frac12 \big)}
\end{gather*}
and
\begin{gather*}
  f^{(2)}(t)  =   \frac{2}{(x^2+x\sqrt{x^2-t}-t)(t+\beta -\frac{1}{2})^3}+\frac{2\left(-1-\frac{x}{2\sqrt{x^2-t}}\right)}{\big(x^2+x\sqrt{x^2-t}-t\big)^2 \big(t+\beta -\frac{1}{2}\big)^2}\nonumber \\
\hphantom{f^{(2)}(t)  =}{}
 +  \frac{2\left(-1-\frac{x}{2\sqrt{x^2-t}}\right)^2}{\big(x^2+x\sqrt{x^2-t}-t\big)^3 \big(t+\beta -\frac{1}{2}\big)} \\
\hphantom{f^{(2)}(t)  =}{}
 +  \frac{x}{4(x^2-t)^{3/2}\big(x^2+x\sqrt{x^2-t}-t\big)^2 \big(t+\beta -\frac{1}{2}\big)}.
\end{gather*}
Based on \eqref{f-def} and the above two formulas, it is readily shown that
\begin{gather}\label{EM-third}
  \frac{1}{2}(f(n)+f(1)) = O\big(x^{-2}\big), \qquad \frac{B_2}{2} \big[f^{(1)}(n) - f^{(1)}(1)\big] =  O\big(x^{-2}\big),
\end{gather}
and
\begin{gather}\label{EM-remainder}
|R_1|\leqslant \big(2-2^{1-2}\big)\frac{|B_2|}{2!}\int_{1}^n \big|f^{(2)}(t)\big| dt =O\big(x^{-2}\big).
\end{gather}
Substituting \eqref{EM-leading}--\eqref{EM-remainder} into \eqref{EM-sum-m=1} yields \eqref{sum-log-n}.
\end{proof}

Combining the results obtained  above, we have the following asymptotic approximation of $\phi_{n}(\sqrt{n+ \beta -1}   z)$ for $z$ bounded away from the transition point 1 and the point origin.

\begin{Proposition}\label{prop-non-uniform}
  Let $N = n+ \beta -1$. Then as $n \to \infty$, we have
    \begin{gather}
    \phi_{n}\big(N^\frac{1}{2} z\big)  =
     \frac{ (4e)^{- \frac{N} {2} } \gamma_n   N^{\frac{n}{2}}  e^{N z^2 + \frac{(\alpha^2-\frac{1}{4}) \log(Nz^2) }{4Nz^2} }   }{z^{\beta - \frac{3}{2}} (z^2-1)^{\frac{1}{4}}} \nonumber \\
\hphantom{\phi_{n}\big(N^\frac{1}{2} z\big)  = }{}
 \times \exp\left\{N \left( \log\big(z + \sqrt{z^2-1} \big)- z \sqrt{z^2-1} \right) \right\}
     \big(1  + O\big(N^{-1}\big)\big) \label{phi-outside}
  \end{gather}
  uniformly for $z$ in any compact subset of ${\mathbb C}{\setminus}[-1,1]$; and
  \begin{gather}
   \frac{\phi_{n}\big(N^\frac{1}{2} z\big)}{(4e)^{- \frac{N} {2} } \gamma_n   N^{\frac{n}{2}}  e^{N z^2 + \frac{(\alpha^2-\frac{1}{4}) \log(Nz^2) }{4Nz^2} } }  \nonumber\\
   \qquad{} =   \frac{2 }{z^{\beta - \frac{3}{2}} (1-z^2)^{\frac{1}{4}}}   \cos \left\{N \left( z \sqrt{1-z^2} - \cos^{-1} z \right) + \frac{\pi}{4} \right\}
    + O\big(N^{-1}\big), \label{phi-inside}
  \end{gather}
  uniformly for $z$ in any small complex neighborhood of $[\delta, 1-\delta]$, where $\delta > 0$ is a fixed small number.
\end{Proposition}

\begin{proof}
  From \eqref{pi-wk}, we have
  \begin{gather*}
    \log \pi_n(x) = \sum_{k=1}^n \log w_k(x).
  \end{gather*}
  The above formula and \eqref{wk-asy} give us, for $|z| > 1$,
  \begin{gather}
    \log \pi_n\big(\sqrt{N}   z \big)   =   n \log\left(\frac{\sqrt{N}}{2}\right) +  \sum_{k=1}^n \log\left(z + \sqrt{z^2 - \frac{2\lambda_k}{N}} \right) \nonumber \\
\hphantom{\log \pi_n\big(\sqrt{N}   z \big)   =}{}
  + \sum_{k=1}^n \log\left[ 1+ \frac{\lambda_k - \lambda_{k-1}}{2[Nz^2 - 2\lambda_k]} \right] + O \big(n^{-1}\big). \label{log-pi-1}
  \end{gather}
  Taking a careful calculation, we f\/ind
  \begin{gather}
  \log\left(z + \sqrt{z^2 - \frac{2\lambda_k}{N}} \right) =  \log\left(z + \sqrt{z^2 - \frac{k}{n}} \right) - \frac{1}{2n} \frac{\beta - \frac{3}{2}}{z^2 - \frac{k}{n} + z \sqrt{ z^2 - \frac{k}{n}} } \nonumber \\
\qquad{}  + \frac{1}{2n} \frac{\alpha^2 - \frac{1}{4}}{\big(k+\beta-\frac{1}{2}\big)  \big[z^2 - \frac{k}{n} + z \sqrt{ z^2 - \frac{k}{n}} \big]} + \frac{\beta-1}{2n} \left(\frac{z }{\sqrt{ z^2 - \frac{k}{n}}} - 1 \right) + O\big(n^{-2}\big)
  \end{gather}
  and
  \begin{gather}
    \log \left[ 1+ \frac{\lambda_k - \lambda_{k-1}}{2[ Nz^2 - 2\lambda_k]} \right] = \frac{1}{4n} \frac{1}{z^2 - \frac{k}{n}} + O\big(n^{-2}\big),
  \end{gather}
  where the error terms $O(n^{-2})$ are uniformly valid for all $k$ in the above two formulas. As $n \to \infty$, by the trapezoidal rule, we have
  \begin{gather}
   \sum_{k=1}^n \log\left(z + \sqrt{z^2 - \frac{k}{n}} \right) \sim n \int_0^1 \log\big(z + \sqrt{z^2 - s } \big) ds  + \frac{1}{2} \log \frac{z + \sqrt{z^2-1}}{2z} \nonumber \\
 \qquad{} = n \left(z^2 - \frac{1}{2} - z \sqrt{z^2 -1} + \log\big(z + \sqrt{z^2-1}\big) \right) + \frac{1}{2} \log \frac{z + \sqrt{z^2-1}}{2z},
  \\
  \sum_{k=1}^n \left[  - \frac{1}{2n} \frac{\beta - \frac{3}{2}}{z^2 - \frac{k}{n} + z \sqrt{ z^2 - \frac{k}{n}} }
   + \frac{\beta-1}{2n} \left(\frac{z }{\sqrt{ z^2 - \frac{k}{n}}} - 1 \right) \right] \nonumber \\
\qquad{} \sim -\frac{\beta - \frac{3}{2}}{2} \int_0^1 \frac{ds}{ z^2-s + z \sqrt{z^2 -s }} + \frac{\beta-1}{2}  \int_0^1 \frac{z ds}{\sqrt{ z^2-s }} - \frac{\beta -1}{2}
    \nonumber \\
\qquad{}  = \left(\beta - \frac{3}{2}\right) \log \frac{z + \sqrt{z^2-1}}{2z} + (\beta - 1) \left( z - z \sqrt{z^2-1} - \frac{1}{2} \right)
  \end{gather}
  and
  \begin{gather} \label{log-pi-last}
    \sum_{k=1}^n \frac{1}{4n} \frac{1}{z^2 - \frac{k}{n}} \sim \frac{1}{4} \int_0^1 \frac{1}{z^2-s} ds = \frac{1}{4} \log\frac{z^2}{z^2-1}.
  \end{gather}
  Combining \eqref{sum-log-n} and the above formulas \eqref{log-pi-1}--\eqref{log-pi-last}, the asymptotic formula \eqref{phi-outside} in the exponential region follows from the relation $\phi_n(x) = \gamma_n \pi_n(x)$.

  The asymptotic formula \eqref{phi-inside} near the oscillatory region can be easily obtained via a mat\-ching method introduced in~\cite{xswang2014,xswang2012}.
  Here, we only sketch the idea of matching principle.
  Denote the leading term on the right-hand side of the asymptotic formula \eqref{phi-outside} by~$\Phi(n,z)$.
  It is readily seen that~$\Phi(n,z)$ is analytic in~${\mathbb C}{\setminus}[-1,1]$.
  Furthermore, $\Phi(n,z)$ can be analytically continued to the interval $(-1,0)\cup(0,1)$ from either above or below (note that $\Phi(n,z)$ has a singularity at~$z=0$). By symmetry, we only need to consider the interval~$(0,1)$.
  Denote
  \[
  \Phi_\pm(n,z):=\lim_{\varepsilon\to0^+}\Phi(n,z\pm i\varepsilon),\qquad 0<z<1.
  \]
  We can write
  \begin{gather*}
    \Phi_\pm(n,z)=\frac{ (4e)^{- \frac{N} {2} } \gamma_n   N^{\frac{n}{2}}  e^{N z^2 + \frac{(\alpha^2-\frac{1}{4}) \log(Nz^2) }{4Nz^2} }   }{z^{\beta - \frac{3}{2}} (1-z^2)^{\frac{1}{4}}e^{\pm i\pi/4}} \times \exp\left\{\pm iN \left( \cos^{-1}z-z \sqrt{1-z^2} \right) \right\}.
  \end{gather*}
  Clearly, both $\Phi_+(n,z)$ and $\Phi_-(n,z)$ are analytic in a small complex neighborhood of $[\delta, 1-\delta]$, where $\delta > 0$ is a f\/ixed small number. It is also noted that as $n\to\infty$, the ratio
  \[
  {\Phi_+(n,z)\over\Phi_-(n,z)}=\exp\left\{ 2iN \left( \cos^{-1}z-z \sqrt{1-z^2} \right)- i\pi/2 \right\}
  \]
  is exponentially large when $z$ is above and close to the interval $[\delta, 1-\delta]$, and exponentially small when $z$ is below and close to the interval $[\delta, 1-\delta]$.
  Therefore, we have
  \[
  \phi_n\big(N^{1\over2}z\big)=\Phi_+(n,z)[1+O(1/n)]+\Phi_-(n,z)[1+O(1/n)]
  \]
  for all $z$ with real part in $[\delta, 1-\delta]$ and small nonzero imaginary part.
  By analytic continuation, the above formula is actually valid in a small complex neighborhood of $[\delta, 1-\delta]$.
  Finally, we obtain \eqref{phi-inside} in view of the following identity:
  \begin{gather*}
    \Phi_+(n,z)+\Phi_-(n,z) =\frac{ (4e)^{- \frac{N} {2} } \gamma_n   N^{\frac{n}{2}}  e^{N z^2 + \frac{(\alpha^2-\frac{1}{4}) \log(Nz^2) }{4Nz^2} }   }{z^{\beta - \frac{3}{2}} (1-z^2)^{\frac{1}{4}}}
    \\
    \hphantom{\Phi_+(n,z)+\Phi_-(n,z) =}{}
    \times 2\cos\left\{\pi/4-N \left( \cos^{-1}z-z \sqrt{1-z^2} \right) \right\}.
  \end{gather*}
  This completes the proof of Proposition \ref{prop-non-uniform}.
\end{proof}

\begin{Remark}
  Note that the recurrence coef\/f\/icients $\lambda_n$ given in \eqref{an-def} are regularly varying at inf\/inity.  Therefore, one may apply the techniques in Geronimo et al.~\cite{Ger:Smi,Ger:Smi:VAs} to obtain the similar non-uniform asymptotic results for complex~$z$ bounded away from the oscillatory interval; see Theorem~3.8 in~\cite{Ger:Smi:VAs} and~\eqref{phi-outside} in this paper. However, because more general cases are considered in~\cite{Ger:Smi:VAs}, their results are only implicitly given in terms of an~integral and a~product of a sequence from~$1$ to~$n$, while our formula~\eqref{phi-outside} is more explicit.
  The proof of~\eqref{phi-outside} given above complements the general results obtained by Geronimo et~al.~\cite{Ger:Smi,Ger:Smi:VAs}. Moveover,
  the asymptotic expansion~\eqref{phi-inside} inside the oscillatory interval and the uniform formulas in Theorem~\ref{main-thm}
  are not given in~\cite{Ger:Smi,Ger:Smi:VAs}.
\end{Remark}

Now we are ready to prove our main results.
\begin{proof}[Proof of Theorem~\ref{main-thm}]
Let $x = N^{\frac{1}{2}}   t$. We have from \eqref{phin-pn} and \eqref{phi-inside}
\begin{gather}
  p_n(x) = \frac{\phi_{n}(x)}{k_n}  \sim \frac{2}{k_n} \frac{ (4e)^{- \frac{N} {2} } \gamma_n   N^{\frac{n}{2}}  e^{N t^2 + \frac{(\alpha^2-\frac{1}{4}) \log(Nt^2) }{4Nt^2} } }{t^{\beta - \frac{3}{2}} (1-t^2)^{\frac{1}{4}}}  \nonumber\\
   \hphantom{p_n(x) = \frac{\phi_{n}(x)}{k_n}  \sim}{}
   \times \cos \left\{N \left( t \sqrt{1-t^2} - \cos^{-1} t \right) + \frac{\pi}{4} \right\}.\label{pn-com1}
\end{gather}
Recalling the asymptotic expansions for Airy functions, as $x \to \infty$
\begin{gather*}
  \Ai(-x) \sim  \frac{1}{\sqrt{\pi} x^{1/4}}  \cos \left( \frac{2}{3} x^{\frac{3}{2}}-\frac{\pi}{4}\right), \qquad  \Bi(-x) \sim  -\frac{1}{\sqrt{\pi} x^{1/4}}  \sin \left( \frac{2}{3} x^{\frac{3}{2}} -\frac{\pi}{4}\right),
\end{gather*}
we have, when $t<1$
\begin{gather}
  p_n(x)  =   C_1(x) P_n(x) + C_2(x) Q_n(x) \sim C_1(x)  \left(\frac{1}{1- t^2 }\right)^{\frac{1}{4}}
   \frac{1}{\sqrt{\pi} }  \cos \left( \frac{2N}{3} (-U(t))^{\frac{3}{2}} -\frac{\pi}{4}\right) \nonumber \\
\hphantom{p_n(x)  =}{} - C_2(x) \left(\frac{1 }{1 - t^2 }\right)^{\frac{1}{4}} \frac{1}{\sqrt{\pi} }  \sin \left( \frac{2N}{3} (-U(t))^{\frac{3}{2}} -\frac{\pi}{4} \right) \nonumber \\
\hphantom{p_n(x)}{} =  C_1(x)  \left(\frac{1 }{1 - t^2 }\right)^{\frac{1}{4}}
   \frac{1}{\sqrt{\pi} }  \cos \left( N\left(\cos^{-1} t -  t \sqrt{1-t^2}\right)  -\frac{\pi}{4}\right) \nonumber \\
\hphantom{p_n(x)  =}{}
 - C_2(x) \left(\frac{1 }{1 - t^2 }\right)^{\frac{1}{4}} \frac{1}{\sqrt{\pi} }   \sin \left( N\left(\cos^{-1} t - t \sqrt{1-t^2}\right)  -\frac{\pi}{4}\right). \label{pn-com2}
\end{gather}
Comparing \eqref{pn-com1} and \eqref{pn-com2}, we have
\begin{gather*}
  C_1(x) = \pi^{\frac{1}{4}} \sqrt{\frac{\Gamma(\beta + \alpha) \Gamma(\beta -\alpha)}{\Gamma(\beta + \frac{1}{2})}} \frac{e^{x^2 + \frac{(\alpha^2-\frac{1}{4}) \log(x^2) }{4x^2} }}{ (2x)^{\beta-\frac{3}{2}}} ,  \qquad C_2(x) = 0.
\end{gather*}
This completes the proof of the theorem.
\end{proof}

\section{Discussions} \label{sec:discuss}

In this paper, we investigated a family of orthogonal polynomials associated with certain coherent states introduced by Ali and Ismail \cite{Ali:Ism}. Using the Wang--Wong dif\/ference equation method, we obtained uniform asymptotic expansions for the polynomials as their degree tends to inf\/inity. It is interesting to note that our asymptotic formula \eqref{main-for} contains a term $x^{(\alpha^2-1/4)/(2x^2)}$, which is generated from the $\frac{\alpha^2-\frac{1}{4}}{n+\beta-\frac{1}{2}}$ term in the expression \eqref{an-new-expression} for the recurrence coef\/f\/icients of three-term recurrence relation. We remark that none of the hypergeometric orthogonal polynomials within the Askey scheme has this type of phenomenon; namely, for any polynomials in the Askey scheme, the $O(1/n)$ term of the series expression for the recurrence coef\/f\/icients of three-term recurrence relation always vanishes. According to our asymptotic results in \eqref{main-for}, this phenomenon also suggests that the weight function for the orthogonal polynomials $\phi_n(x)$ may contain a factor $|x|^{(\frac{1}{4}-\alpha^2)/x^2}$. Furthermore, we conjecture that the weight function will violate the Szeg\H{o} condition at the origin due to this factor $|x|^{(\frac{1}{4}-\alpha^2)/x^2}$. Note that such kind of factor is also dif\/ferent from the known orthogonal polynomials in the non-Szeg\H{o} class, such as the Pollaczek polynomials in \cite[p.~393]{Sze}.

Several questions remain open.
\begin{enumerate}\itemsep=0pt
\item[(i)] One may try to f\/ind an explicit formula for the weight function of $\phi_n(x)$, then discover an intrinsic relation between the two sets of polynomials $\phi_n(x)$ and $\psi_n(x)$.

\item[(ii)] One can see that our formula \eqref{main-for} has singularities near the origin if $\alpha\neq1/2$. New ideas are needed to  derive uniform asymptotic formulas for $\phi_n(N^{1/2}t)$ when $t\in[-\delta,\delta]$ and $\alpha\neq1/2$.

\item[(iii)] As the weight function for $\psi_n(x)$ is given explicitly in~\eqref{psi-weight}, it is possible to adopt the Deift--Zhou nonlinear steepest descent method for Riemann--Hilbert problems to f\/ind asymptotic formulas for~$\psi_n(x)$.
\end{enumerate}

\subsection*{Acknowledgements}

We are grateful to Professor Mourad Ismail for his helpful comments and discussions.
We would also like to thank the anonymous referees for their constructive suggestions which lead to a~signif\/icant improvement of the manuscript.

DD and WH were partially supported by grants from the Research Grants Council of the Hong Kong Special Administrative Region, China (Project No.~CityU 101411, CityU 11300814). In addition, DD was partially supported by a grant from the City University of Hong Kong (Project No.~7004065). XSW was partially supported by two grants (Start-up Fund and Summer Research Grant) from Southeast Missouri State University.

\pdfbookmark[1]{References}{ref}
\LastPageEnding

\end{document}